\documentclass[11pt,reqno]{article}

\usepackage[usenames]{color}
\usepackage{amssymb}
\usepackage{amscd}
\usepackage{graphicx} 
\usepackage[all]{xy}

\usepackage{amsmath,amssymb}
\usepackage{amsthm}
\usepackage{amsfonts}
\usepackage{latexsym}
\usepackage{url}

\begin{document}

\begin{center}
{\LARGE
Spheres of small diameter with long sweep-outs\\
}
\vskip 0.1 in
\large
Yevgeny Liokumovich
\footnote{The author's research is supported by NSERC CGS scholarship}
\\
\end{center}
\vskip 0.1in

\theoremstyle{plain}
\newtheorem{theorem}{Theorem}
\newtheorem{lemma}[theorem]{Lemma}
\newtheorem{remark}[theorem]{Remark}

\newcommand{\N}{\mathbb N}

\begin{abstract}

We prove the absence of a universal diameter bound on lengths of curves in a sweep-out of a Riemannian 2-sphere. If such bound existed it would yield a simple proof of existence of short geodesic segments and closed geodesics on a sphere of small diameter.

\end{abstract}

\section{Introduction} 
By a sweep-out of a Riemannian 2-sphere $M=(S^2,g)$ we mean a non-contractible loop $\gamma_t$ in $(\Lambda M,\Lambda ^0 M)$, where $\Lambda M$ denotes the space of free loops on $M$ and $\Lambda ^0 M$ denotes the space of constant loops. In other words, $\gamma_t$ is a 1-parameter family of closed curves, starting and ending at a point and inducing a non-zero degree map $f:S^2 \rightarrow M$. If there exists a sweep-out by loops of length $\leq L$ then the standard minimax argument implies that there exists a non-trivial closed geodesic on $M$ of length $\leq L$ ([B]).

Moreover, if there exits a sweep-out of $M$ by loops with 2 fixed points and of length $\leq L$ then for any two points $a,b \in M$ there exists $k$ distinct geodesic segments from $a$ to $b$ of length $\leq 2kL + 2diam(M)$ (cf. [NR1]). If we could take $L \leq C diam(M)$ we would obtain a linear in $k$ bound on the length of $k$-th shortest geodesic segment in terms of the diameter. 

Yet A.Nabutovsky and R.Rotman observe (see [N], [NR1], [NR2]) that example of S.Frankel and M.Katz [FK] (see below) suggests that for any $C$ there exists a sphere with no sweep-out obeying this inequality. The author learnt from Regina Rotman that a similar conjecture was independently made by S.Sabourau. In this note we prove this conjecture.

\section{Main result} 

Let $\gamma: I \times S^1 \rightarrow M$ be a 1-parameter family of free loops (we write $\gamma_t(s)$ for $\gamma(t,s)$), such that $\gamma_0(s)$ and $\gamma_1(s)$ are constant loops. $\gamma_t$ induces a map $f:S^2 \rightarrow M$, such that $\gamma(t,s)=f \circ p(t,s)$, where $p: I \times S^1 \rightarrow S^2$ is the suspension map that collapses $\{0\} \times S^1$ to the South pole and $\{1\} \times S^1$ to the North pole of $S^2$. If $deg(f) \neq 0$ we call $\gamma_t$ a sweep-out.

\begin{theorem}
For any $C>0$ there exists a Riemannian 2-sphere $M$ of diameter $\leq 1$, such that for any sweep-out $\gamma_t$ of $M$ there is a loop $\gamma_{t_0}$ of length $\geq C$.
\end{theorem}

F.Balacheff and S.Sabourau in [BS] defined a \emph{diastole} of $M$ as 

$$dias(M):= inf_{(\gamma_t)} sup _{0 \leq t \leq 1} length(\gamma_t) $$

where $(\gamma_t)$ runs over the families of loops inducing $f:S^2 \rightarrow M$ of degree $\pm 1$. In [S, Remark 4.10] S.Sabourau constructs a sequence $M_n$ of Reimannain two-spheres such that

$$lim_{n \rightarrow \infty} \frac{\sqrt{Area(M_n)}}{dias(M_n)} = 0.$$ 

Theorem 1 implies the analogous result with $daim (M_n)$ in place of $\sqrt{Area(M_n)}$.

F.Balacheff and S.Sabourau prove ([BS]) that if 1-parameter families of loops in the definition of the diastole are replaced with 1-parameter families of one-cycles then every Riemannian 2-sphere satisfies $dias(M)\leq C \sqrt{Area(M)}$ for a universal constant $C$.

\begin{proof}
We use the example of S.Frankel and M.Katz [FK]. For any natural number $N$ they embed a binary tree $T$ of height $N$ in a 2-dimensional disc $D$ and define a Riemannian metric on $D$, such that the distance between any two non-adjacent edges of $T$ $dist(e_i,e_j) \geq 1/2$, but the diameter $diam(D) \leq 1$. They prove that for every homotopy of closed curves $\gamma_t$ with $\gamma_0=\partial D$ and $\gamma_1 = \{*\}$ there is an intermediate curve $\gamma_{t_0}$ that intersects at least $O(N/log N)$ edges of $T$ and hence must be at least $O(N/log N)$ long.

Let $M=(S^2,g)$ be a sphere of diameter less than 1 containing the disc of Frankel and Katz with an embedded binary tree $T$ of height $N$. 
Consider a sweep-out $\gamma_t$ and let $f:S^2 \rightarrow M$ be the induced map from the suspension of $S^1$ to $M$ ($\gamma_t(s)=f \circ p(t,s)$).

Suppose that at time $t_0 \in [0,1]$ $\gamma_{t_0}$ does not pass through a vertex $v \in T$. Let $E$ denote a connected component of $M \setminus \gamma_{t_0}$ that contains $v$. Identify all points of $M \setminus E$ to a new point $x$ to obtain a quotient space $E \cup \{ x \}$ homeomorphic to  $S^2$. The loop $p(t_0,s)$ divides $S^2$ into two connected components, the south component $S_{t_0}=\{p(t,s): t< t_0\}$ and the north component $N_{t_0}=\{p(t,s): t> t_0\}$. Under the composition of $f$ with the quotient map $q_1:M \rightarrow E \cup \{ x \}$ the loop $p(t_0,s)$ is mapped to the point $x$. If we collapse all points of $p(t_0,s)$ to a point $y$ we obtain a map $f_{t_0,v}: S_{t_0} \cup \{ y \} \cup N_{t_0} \rightarrow E \cup \{ q \}$, where the quotient space $S_{t_0} \cup \{ y \} \cup N_{t_0}$ is homeomorphic to the wedge sum $S^2 \vee S^2$. We have the following commutative diagram

\begin{center}
\
\xymatrix{
S^2 \ar[d]_{q_2} \ar[r]^f & M \ar[d]^{q_1} \\
S_{t_0} \cup \{ y \} \cup N_{t_0} \ar[r]^{f_{t_0,v}} & E \cup \{ x \}
}
\
\end{center}

Let $f^S _{t_0,v}$ denote the restriction of $f_{t_0,v}$ to the south sphere $S_{t_0} \cup \{ y \}$  
and $f^N _{t_0,v}$ denote the restriction to the north sphere $N_{t_0} \cup \{ y \}$. 
Observe that from the induced commutative diagram for the second homology groups we have 
$deg(f)=deg(f^S _{t_0,v})+deg(f^N _{t_0,v})$. Indeed, the map 
$q_{2 \sharp} : H_2(S^2) \rightarrow H_2(S^2 \vee S^2)$ 
sends a generator $1$ to an element $(1,1) \in \mathbb{Z} \times \mathbb{Z}$ and 
$(f_{t_0,v})_{\sharp} (a,b)=(f^S _{t_0,v})_\sharp (a) + (f^N _{t_0,v})_\sharp (b)$, 
while $q_{1 \sharp}$ is an isomorphism.

Let $A \subset [0,1]$ denote the set of all $t$ such that $gamma_t$ does not pass through $v$. We define a function $d_v:A \rightarrow \N$ (degree of $v$ at time $t$) by 

\begin{equation}
d_v(t)= deg(f^N _{t,v})
\end{equation}

We need two simple facts about $d_v(t)$:

Observation 1. If $\gamma_t$ does not intersect $v$ for $t \in [t_1,t_2]$ then $d_v(t)$ is constant on $[t_1,t_2]$

\emph{Proof:}
Choose a small disc $D$ around $v$, s.t. $D \cap f(p([t_1,t_2])$ is empty and define a quotient map $q$ collapsing $M \setminus D$ to a point. For each $t \in [t_1,t_2]$ we can define a map between spheres $q':N_t \cup \{a \} \rightarrow N_{t_2} \cup \{ b \}$ sending $N_t \setminus N_{t_2}$ to $b$. Then $q \circ f^N _{t,v}= q \circ f^N _{t_2,v} \circ q'$. As $deg(q)=deg(q')=1$ we have $deg(f^N _{t,v})=deg( f^N _{t_2,v} )$

Observation 2. Let $v_1$, $v_2$ be two vertices of $T$ connected by an edge $e$. If $\gamma_t$ does not intersect $e$ then $d_{v_1}(t)=d_{v_2}(t)$

\emph{Proof:}
Since $\gamma_t$ does not intersect $e$ we have that $v_1$ and $v_2$ belong to the same connected component of $M \setminus \gamma_t$. Hence, $f^N _{t,v_1}=f^N _{t,v_2}$.

Without loss of generality we may assume that the images of North and South poles under $f$ are not vertices of $T$. Now we observe that $d_v(0)=deg(f) \neq 0$ and $d_v(1)=0$ for all vertices $v$.

The rest of the proof proceeds as in [FK]. Let $V$ be the set of vertices of $T$ and $K(t)=\# \{v \in V : d_v=0\}$. We may perturb the homotopy slightly, so that $\gamma_t$ passes through no more than one vertex for each $t$. Observation 1 implies that as $t$ varies between $0$ and $1$ $K(t)$ will attain every value between $0$ and $2^N-1$ (recall that $N$ is the height of $T$). 

Consider what values $K(t)$ can attain if $\gamma_t$ intersects only one edge of $T$. Let $v_1$ and $v_2$ be two vertices of an edge $e$ at distances (in the standard metric on the tree) $i$ and $i+1$, respectively, from the root. By Observation 2 we have the following possibilities for the value of $K(t)$: the number of vertices in the connected component of $M \setminus \gamma_t$ that contains $v_1$ ($2^N-2^{N-i-1}$), the number of vertices in the connected component that contains $v_2$ ($2^{N-i-1}-1$) or one of the min or max values ($0$ and $2^N-1$). Since this is true for every d $i$ $K(t)$ can attain at most $2(N-1)+2$ possible different values if $\gamma_t$ intersects only one edge of $T$. Similarly, if have exactly $j$ intersections then there are at most $(2N)^j$ choices for the value of $K(t)$. If throughout the homotopy $\gamma_t$ intersects at most $k$ edges then $K(t)$ attains no more than $\sum _{j=1} ^k (2N)^j \leq (2N)^{k+1}$ distinct values. Since all possible values are attained, we have $(2N)^{k+1} \geq 2^N-1$ and hence $k \geq O(N /log N)$. 

Since the distance between any two non-adjacent edges is greater than $1/2$ we have that for some $t$ $\gamma_t$ will be at least $O(N /log N)$ long.
\end{proof}

\section{Acknowledgements}
The author gratefully acknowledges the support by Natural Sciences and Engineering Research Council (NSERC) CGS scholarship. The author would like to thank Professors Alexander Nabutovsky and Regina Rotman for explaining the question to him and for many very helpful discussions.

\begin{scshape}
Department of Mathematics, University of Toronto, Toronto, Ontario M5S 2E4, Canada 
\end{scshape}

\emph{E-mail address}: \url{e.liokumovich@utoronto.ca}

\end{document}